\documentclass[11pt]{article}

\usepackage{amsmath, amssymb}
\usepackage{graphicx}
\newcommand {\R} {\mathbb{R}}
\newcommand {\Z} {\mathbb{Z}}
\newcommand {\N} {\mathbb{N}}

\newtheorem{Theo}{Theorem}[section]
\newtheorem{Prop}{Proposition}[section]
\newtheorem{Cor}{Corollary}[section]
\newtheorem{Def}{Definition}[section]

\begin{document}
\parindent=1cm

\centerline{\Large Rectangles, integer vectors and hyperplanes of the hypercube}\vskip 8mm

\centerline{\bf E. Gioan , I. P. Silva}\vskip 1cm


{\parindent=0cm
{\small \bf Abstract.}{\small We introduce a family of nonnegative integer vectors - primitive vectors - defining hyperplanes of the real affine cube over  $C^n:=\{-1,1\}^n$  and study their properties with respect to the rectangles of the cube. As a consequence we give a short proof that, for small dimensions ($n\leq 7$),  the real  affine cube can be recovered from its signed rectangles and its signed cocircuits complementary of its facets and skew-facets.}  } \vskip 1.5cm

\section {Introduction}

We consider  $C^n:=\{-1,1\}^n$. $\mathcal Aff(C^n)$  denotes the oriented matroid of the affine dependencies of  $C^n$  over  $\R$, {\it the real affine cube}.

The oriented matroid  $\mathcal Aff(C^n)$  can be defined in several equivalent ways. We recall its definition in terms of the family of signed cocircuits. 

A  signed cocircuit of  $\mathcal Aff(C^n)$  is an ordered pair  $X=(X^+,X^-)$  satisfying the following two conditions: 
{\parindent=0cm

1) $H:= C^n\setminus (X^+\cup X^-)$ is the set of all elements of  $C^n$  spanning the same affine hyperplane   $H$ of  $\R^n$. 

2)  If  the hyperplane   $H:= C^n\setminus (X^+\cup X^-)$  is defined by a linear equation  $H:{\bf x.h}=b$  $({\bf h},b)\in \Z^{n+1}$  then  either $X^+:=\{ {\bf v}\in C^n : {\bf v.h}>b\}$  and  $X^-:=\{ {\bf v}\in C^n : {\bf v.h}<b\}$   or  $X^+:=\{{\bf  v}\in C^n : {\bf v.h}<b\}$  and  $X^-:=\{{\bf  v}\in C^n : {\bf v.h}>b\}$.  }

To every hyperplane  $H:{\bf x.h}=b$  of the matroid  $\mathcal  Aff(C^n)$ is associated the pair of opposite signed cocircuits,   $X=(X^+, X^-)$  and  $-X=(X^-,X^+)$. This pair encodes the  $2$-partition of  $C^n\setminus H= X^+\uplus X^-$  into the points lying on each one of the open half-spaces  of  $\R^n$  defined by the affine hyperplane spanned by H.   

The {\it oriented matroid   $\mathcal Aff(C^n)$}  is defined by the collection  of all its signed cocircuits. 

We suggest that the reader consults  $\cite {OM}$  as general reference on oriented matroids and   $\cite{Ox}$   as general reference on matroids.

The {\it (signed) rectangles}  of  $C^n$  are the shortest (signed) circuits of  $\mathcal Aff(C^n)$. They are the signed sets of the  form    $\pm R$  with  $R=(\{{\bf u,v}\}, \{{\bf u',v'}\})={\bf u^+v^+u^{'-}v^{'-}}$,  where  $({\bf u,u',v,v'})$ are the vertices of a rectangle  of  $\R^n$  with diagonals  $\bf uv$,  $\bf u'v'$.  We denote by  $\mathcal R$  the family of the signed rectangles the real affine cube.

The  {\it Facets} of  $C^n$  are the subsets:  $H_{i^+}:= \{ {\bf v}\in C^n :  v_i=1 \}$   and   $H_{i^-}:= \{ {\bf v}\in C^n :  v_i=-1 \},\ i=1,\ldots, n.$  The {\it skew-facets} of  $C^n$  are the subsets:  $H_{ij^+}:= \{ {\bf v}\in C^n :  v_i=v_j \}$   and  $H_{ij-}:= \{ {\bf v}\in C^n :  v_i=-v_j \},  \ 1\leq i<j\leq n.$. They are hyperplanes of the real affine cube as well as its shortest cocircuits. 

The signed cocircuits of  $\mathcal Aff(C^n)$   complementary of the facets and of the skew facets are  denoted  respectively:  $X_{i+}=( H_{i-},\emptyset)$,  $X_{i-}=( H_{i+},\emptyset)$   and  $X_{ij^+}=(H_{i^+}\cap H_{j^-}, H_{i^-}\cap H_{j^+})$,  $X_{ij^-}=(H_{i^+}\cap H_{j^+}, H_{i^-}\cap H_{j^-})$.

We denote   $\mathcal F$  the family of  the signed cocircuits complementary of the facets of the skew facets of the real cube  $\mathcal Aff(C^n)$.

\begin{Def} {  An oriented cube ( canonically oriented cube of (\cite{dS2})) is an oriented matroid over  $C^n$  containing as signed circuits the signed rectangles of $\mathcal R$  and, as signed cocircuits, the signed cocircuits of  $\mathcal F$.}\end{Def}

In this note we give a short proof of the following theorem (Theorem 3.1) whose proof is mentioned in (\cite {dS2}).\vskip 4mm

{\parindent=0cm
{\bf Theorem} {\it  For $n\leq 7$  the real affine cube   $\mathcal Aff(C^n)$  is the unique oriented cube.} }\vskip 4mm

This is a small step towards an answer to the next Question 1. A positive answer to this question with the results of \cite{dS2} would imply the existence of a purely combinatorial characterization of affine/linear dependencies of  $\pm 1$  vectors over the reals:\vskip 2mm

{\parindent=0cm
{\bf Question.}(\cite {dS2}){\it Is the real cube the unique oriented cube?}}\vskip 2mm

Our proof of the theorem uses the encoding of  $\mathcal Aff(C^n)$  in terms of the family  $\mathcal H_ n$  of non-negative integer vectors of  $\N_0^{n+1}$  that defines its hyperplanes and signed cocircuits up to symmetries of the non-oriented matroid of the real affine cube  (see \cite{dS1}). This is briefly recalled in section 3, where the complete proof is presented. 

In the next section 2 we define {\it primitive vectors of  $C^n$}  and from there  a recursive family of hyperplanes - {\it primitive hyperplanes of 
$\mathcal Aff(C^n)$ } whose behaviour with respect to the net of signed rectangles implies that the corresponding signed cocircuits of the real cube must be signed cocircuits of every oriented cube. 

We point out that as a direct consequence of the main theorem we obtain (for   $n\leq 7$) a new proof of the following conjecture of M. Las Vergnas also open for  $n>7$:\vskip 2mm

{\parindent=0cm
{\bf Las Vergnas cube conjecture (\cite{LV}):} {\it $\mathcal Aff(C^n)$  is the unique orientation of the (non-oriented) real affine cube.}}\vskip 2mm

This conjecture was verified computationally for  $n\leq 7$, by J. Bokowski etal \cite{BG}  .

\section { Nonnegative integer vectors, rectangles and hyperplanes of oriented cubes}

{\bf 1.1. Parallel Strata and rectangles}\vskip 2mm

We are interested in nonnegative integer vectors  ${\bf h}\in \N_0^n$  and how they stratify the vertices of the cube  $C^n\in \R^n$ into parallel  levels (strata)  orthogonal to  $\bf h$. 

{\parindent=0cm
\begin{Def}{\rm ( $\bf h$-levels of  $C^n$)}

Given  ${\bf h}\in \N_0^n$  we consider  $|{\bf h}|:=\sum _{i=1}^n h_i$.

For every  $a\in \{0,\ldots, |{\bf h}|\}$  the  $a$-level of $\bf h$, denoted  $S_a({\bf h})$  or simply  $S_a$,  is the set of vertices of the cube   $C^n$  defined as:
$$S_a=S_a({\bf h}):= \{ {\bf v}\in C^n: \ {\bf h. v}= |{\bf h}|- 2a \}$$
Clearly   $S_{|h|-a}= - S_a$, in particular  $S_0=\{{\bf 1}\}$  and  $S_{|{\bf h}|}=\{{\bf -1}\}$. 

Eventually  $S_a =\emptyset$.  \end{Def}

\begin{Def}{\rm ( Realizable  $\bf h$ - rectangles)}

An $\bf h$ - rectangle is a sequence of four natural numbers  $r=(a\leq b \leq c\leq d)$  such that  $d=b+c-a$  and  $0\leq a$  and  $d\leq |{ \bf h}|$.  We distinguish between  $1,\ 2,\ 3$ and $4$- rectangles according to the number of different $\bf h$-levels they cross.

A  $1$-rectangle is of the form  $r_1= (a=a=a=a)$, a  $2$-rectangle is of the form  $r_2= (a=a<b=b)$, $3$-rectangle is of the form  $r_3= (a<b=b<c)$  and a  $4$-rectangle is of the form  $r_4= (a<b<c<d)$,

An  $\bf h$-rectangle  $r=(a\leq b\leq c\leq d)$  is  realizable if there is a signed geometric rectangle, $R={\bf v}_a^+{\bf v}_b^-{\bf v}_c^-{\bf v}_d^+$, of  $\mathcal Aff(C^n)$  such that  ${\bf v}_a\in S_ a({\bf h})$, ${\bf v}_b\in S_ b({\bf h})$, ${\bf v}_c\in S_ c({\bf h})$  and  ${\bf v}_d\in S_ d({\bf h})$. 
\end{Def}

\begin{figure}[h] 
\begin{center}
\includegraphics[scale=0.60]{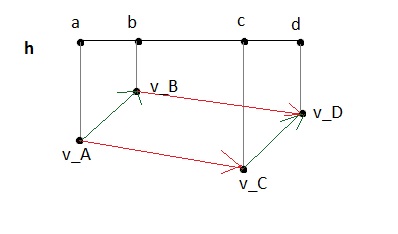}
\end{center}
\vskip - 2mm
\caption{ }
\end{figure}

\begin{Def}{\rm (Embedding of levels)}
Given two   ${\bf h}-levels$,  $S_a$, $S_b$  we say that  $S_a$  is embedded in  $S_b$, written  $S_a \hookrightarrow S_b$   , if there is a spanning tree  $\Gamma_a$  of the complete graph  $K(S_a)$  such that for every edge  $\{{\bf u,u'}\}$  of $\Gamma_a\subseteq S_a$ there is a pair of vertices  ${\bf v,v'}\in S_b$  such that  $R={\bf u^+u^{'-}v^{'-}v^+}$  is a signed rectangle of the cube.
\end{Def}

{\bf Remark.}  Notice that  if  $S_a$  is embedded in  $S_b$, $S_a \hookrightarrow S_b$, then the affine span of  $S_a$  must be parallel to the affine span of  $S_b$. The notion of embedding of levels is from this point view a combinatorial version of parallelism in $\R^n$.\vskip 5mm

{\bf Notation.} We usually identify an element  ${\bf v}\in C^n$  with the subset  $\alpha\subseteq [n]$  of its negative entries, more precisely with the sequence of the elements of  $\alpha$ written by increasing order. Example:  $(1,1,1,1)\in C^4 \equiv \emptyset$,  $(1,-1,-1,1)\in C^4 \equiv 23$. 

Let  ${\bf h}\in \N_0^{n-1}$  and  ${\bf g}=({\bf  h},g)\in \N_0^{n}$. The  ${\bf g}$-levels   and  $\bf h$-levels are related in the following way, for every  $a\in \N_0$,  $0\leq a\leq |{\bf g}|$ :
$$S_a({\bf g})=  S_{a}({\bf h})n^+\cup   S_{a-g}({\bf h})n^-$$
where  $S_{a}({\bf h})n^+ = \{({\bf v},1)\in C^n: \ {\bf v}\in S_a({\bf h})\}$  and  $S_{a-g}({\bf h})n^- = \{({\bf v},-1)\in C^n:  {\bf v}\in S_{a-g}({\bf h})\}$. \vskip 5mm

{\parindent=0cm
{\bf Examples} 1)  The vector  $\bf h=(1,1,1,2)$  has  $6$  levels:  $S_0=\emptyset$,  $S_1=\{1,2,3\}$, $S_2=\{12,13,23,4\}$, $S_3=\{ 14,24,34,123\}$,  $S_4=\{ 124,134,234\}$,  and  $S_5=\{ 1234\}$. Observe that for every  $a\not=2$  one has  $S_a\hookrightarrow S_2$, as represented in Figure 2.  Note also that  $S_2\hookrightarrow S_3$  but  $S_2\not\hookrightarrow S_4$. 
\vskip 2mm
2) The vector  $\bf h=(1,2,2,3)$  has  $9$  levels. No level  $S_b=S_b({\bf h})$  satisfies the property   $\forall a\not=2, \  S_a\hookrightarrow S_b$.

\begin{figure}[h] 
\begin{center}
\includegraphics[scale=0.60]{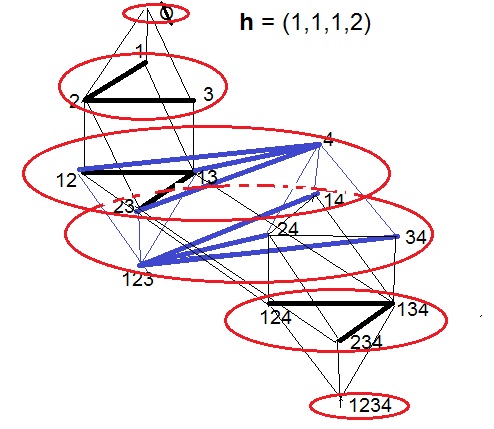}
\end{center}
\vskip - 2mm
\caption{ The (1,1,1,2)-levels and embeddings of the levels  $S_a,\ a\not=2$  in  $S_2$. }
\end{figure}

\begin{Def}{\rm (primitive vectors)}

A vector  ${\bf h}\in \N_0^n$  is called {\it a primitive vector} if all its 3 and 4- numerical rectangles are realizable. 
\end{Def}

{\bf Examples.} The vectors  $(1,1,1)$,  $(1,1,2)$,  $(0,1,1)$,  are primitive vectors of  $C^3$.\vskip 5mm

The next Proposition gives, in particular, a recursive construction of primitive vectors of  $\R^n$. }

{\parindent=0cm
\begin{Prop} {\rm (Properties of primitive vectors)}

Let  ${\bf h}=(h_1,\ldots,h_n)$  be a primitive vector of  $C^n$. Then:

(i) For every  $0\leq a\leq |{\bf h}|$, the  $a$-level  $S_a$  of  $\bf h$  is nonempty , moreover  for  $a\not= 0,|{\bf h}|$, every level  $S_a$  contains at least two elements.

(ii) For every  $g\in \N$,  $0\leq g\leq  \frac{|{\bf h}|}{2}+1$,  ${\bf g}= ({\bf h},g)$  is a primitive vector of  $C^{n+1}$.
\end{Prop}

{\bf Proof.}(i) Imediate from the definition, since every 3-rectangle must be realizable.

(ii) We have  $S_a({\bf g})=  S_{a}({\bf h})(n+1)^+\cup   S_{a-g}({\bf h})(n+1)^-$.}

Let  $r= a<b\leq c< d$  be  a  $3$  or  $4$- rectangle of  ${\bf g}=({\bf h}, g)$. If  $d< |{\bf h}|$  or  $a> g$  the rectangle  $r$ is certainly  $\bf h$- realizable. In the first case with elements of  $S_{a}({\bf h})(n+1)^+$, in the second case with elements of  $S_{a-g}({\bf h})(n+1)^-$ . 

If  $a<g$  and  $d> |{\bf h}|$, then, since  $a+d=b+c >|{\bf h}|$ we can guarantee that  $c\geq g$. Consider the rectangle   $r'=( a<b; c-g<d-g)$. If  $r'$  is a 2-rectangle then  $a=c-g$  and $b=d-g$. Take  ${\bf v}_a\in S_{a}({\bf h})(n+1)^+$  and  ${\bf v}_b\in  S_{b}({\bf h})(n+1)^+$. Then  ${\bf v}_c=({\bf v}_a,-1)\in  S_c({\bf g})$  and  ${\bf v}_d=({\bf v}_b, -1)\in  S_d({\bf g})$  and   $R'=(({\bf v}_a,1), ({\bf v}_b,1); ({\bf v}_c,1), ({\bf v}_d, 1))$. Similarly, if   $r'$  is a  3 or 4  $\bf h$-rectangle then it is realizable in  $S_{a}({\bf h})(n+1)^+$.  Let  $R'=(({\bf v}_a,1), ({\bf v}_b,1); ({\bf v}_{c-g},1), ({\bf v}_{d-g}, 1))$   be a realization of  $r'$  then  $R=(({\bf v}_a,1), ({\bf v}_b,1); ({\bf v}_{c-g},-1), ({\bf v}_{d-g},- 1))$  must be a realization of  $r$  in  $C^n$.

{\parindent=0cm
\begin{Theo} {\rm (primitive vectors and oriented cubes)}

Let $\mathcal M= \mathcal M(C^n)$  be an oriented cube and  $cl: 2^{C^n}\longrightarrow  2^{C^n}$  its closure operator.

Assume that there is a primitive vector  $\bf h$  of  $C^n$  such that some  $\bf h$-level  $S_b$  has the following property:
$$\forall  a\not=b \ \ S_a \hookrightarrow S_b\leqno{(E)}$$ 
Then either  

(i)  $cl(S_b)=C^n$,  or  

(ii) $cl(S_b)=S_b$ is a hyperplane of  $\mathcal M$  whose complement is the support of the pair  $\pm X_b$  of signed cocircuits of  $\mathcal M$  defined by  $X_b= ( \cup_{a<b} S_a,\ \cup_{a>b} S_a)$.
\end{Theo}

{\bf Proof.} First notice that,  by definition of embedding, the hypothesis that $S_b$  satisfies condition  $(E)$  implies that if one element  ${\bf v}\in S_a,\ a\not=b$  belongs to  $cl(S_b)$  then  $S_a\subset cl(S_b)$. On the other hand if some element  ${\bf v}\in S_a,\ a\not=b$  is not  $cl(S_b)$  then, by orthogonality with the signed rectangles of the embedding,  $R=({\bf u}_a^+, {\bf u}^{'-}_a,{\bf v}^{'-}_b , {\bf v}_b^+)$, all the elements of  $S_a$  must have the same sign in any signed covector  $V=(V^+,V^-)$, complementary of  $cl(S_b)$, i.e. one must have either  $S_a\subseteq V^+$  or   $S_a\subseteq V^-$. With these remarks in mind we now prove the theorem.

{\bf Claim 1.} {\it if there is  ${\bf v}\in cl (S_b)$  such that  ${\bf v}\in S_{b-1}\cup S_{b+1}$  then  $cl(S_b)=C^n$.} }

We may assume  $n>1$  and therefore  $b\not= 0,|h|$. 

Assume that there is  ${\bf v}\in S_{b-1}\cap  cl(S_b)$.  In this case, by the above remarks  $S_{b-1}\subset cl(S_b)$. And because the $3$-rectangle  $r=(b-1<b=b<b+1)$  is realizable, orthogonality of the signed covectors complementary of  $cl(S_b)$  with any geometric realization   $R={\bf u}_{b-1}^+ {\bf u}_b^-{\bf v}_b^- {\bf v}_{b+1}^+$  of  $r$, implies that  ${\bf v}_{b+1}$ must be in  $cl(S_b)$  and therefore   $S_{b+1}\subseteq cl(S_b)$. 

Similarly if  ${\bf v}\in S_{b+1}\cap  cl(S_b)$  orthogonality with the same geometric realization  $R$  of the 3-rectangle  $r=(b-1<b=b<b+1)$ implies that  $S_{b+1}\cup  S_{b-1}\subseteq cl(S_b)$. So, if some  ${\bf v}\in  S_{b-1}\cup S_{b+1}$  in the closure of  $S_b$  then  both levels  $S_{b-1}$  and  $S_{b+1}$  are contained in  $cl(S_b)$.

Now, once  $S_{b+1}\cup  S_{b-1}\subseteq cl(S_b)$, considering the sequence of geometric realizations of the $3$-rectangles   $r_b=(b<b+1=b+1<b+2)$,  $r_{b+1}= (b+1<b+2=b+2 < b+3)$, $ \ldots $, $r_{n-2}=(n-2<n-1=n-1<n)$ we conclude, successively, that  $S_{b+2}\subset cl(S_b)\ldots $  $S_n\subset cl (S_b)$.  Similarly using a sequence of 3-rectangles  $r_{b-2}=(b-2<b-1=b-1<b)$, $\ldots$, $r_{2}=(2<1=1<0)$
we also deduce that  $S_{b-2},\ldots, S_0 \subset cl(S_b)$  and therefore that  $cl(S_b)=C^n$.

{\parindent=0cm
{\bf Claim 2.}  {\it If  $(S_{b-1}\cup S_{b+1}) \cap  cl(S_b)=\emptyset$  then  $cl(S_b)= S_b$.}}   

In this case we use the fact that  $\mathcal M$ is a cube canonically oriented  and prove that the only way of signing the complement  $C^n\setminus S_b$,  of  $S_b$,  orthogonally to the rectangles of  $\mathcal M$  is  $\pm X_b$  where  $X_b= ( \cup_{a<b} S_a,\ \cup_{a>b} S_a)$. This proves simultaneously that  $cl(S_b)=S_b$  and that  $X_b$  must be a cocircuit of  $\mathcal M$.

If  $(S_{b-1}\cup S_{b+1}) \cap  cl(S_b)=\emptyset$  then  $cl(S_b)$  is a flat of  $\mathcal M$  whose complement  is the support of (signed) covetors of  $\mathcal M$. Using a realization  $R=({\bf u}_{b-1}^+, {\bf u}_b^-, {\bf v}_b^-, {\bf v}_{b+1}^+)$  of the 3-rectangle  $r=(b-1<b=b<b+1)$  we conclude that $S_{b-1}$  and  $S_{b+1}$  must have different signs in any covector  $V=(V^+,V^-)$  complementary of  $cl(S_b)$. We may assume w.l.o.g.  that  $S_{b-1}\subseteq V^+$  and  that  $S_{b+1}\subseteq V^-$.

Consider the sequence of $3$-rectangles with the "upper" vertex in $S_b$  or  $S_{b-1}$  defined by:  $r_{b+2}= (b<b+1=b+1<b+2)$,  $s_{b+3}=(b-1< b+1=b+1 < b+3)$,  $r_{b+4}= (b<b+2=b+2<b+4)$,  $s_{b+5}=(b-1< b+2=b+2 < b+5)$, $\ldots$ . Orthogonality with geometric realizations of these rectangles implies successively that  $S_{b+2}, S_{b+3}\ldots , S_{|{\bf h}|} \subseteq V ^-$ 

Similarly orthogonality with realizations of the sequence of $3$-rectangles with "lower" vertex in $S_b$  or  $S_{b+1}$  defined by:  $r_{b-2}= (b-2<b-1=b-1<b)$,  $s_{b-3}=(b-3< b-1=b-1 < b+1)$,  $r_{b-4}= (b-4<b-2=b-2<b)$,  $s_{b-5}=(b-5< b-2=b-2 < b+1)$, $\ldots$ . leads us to conclude successively that $S_{b-2}, S_{b-3}\ldots , S_{\emptyset} \subseteq V ^+$. That means there is a unique signature of the complement of  $S_b$  orthogonal to the signed rectangles of  $\mathcal M$ that signature is  defined by  $X_b$. The unicity of the signature and the fact that its support is the complement of  $S_b$  imply that  $X_b$  is a signed cocircuit of  $\mathcal M$  and also that  $cl(S_b)=S_b$  is a hyperplane of  $\mathcal M$.  $\ \square$ \vskip  1cm

The next Proposition gives a recursive way of constructing hyperplanes and the corresponding signed cocircuits of any oriented cube.\vfill\eject

{\parindent=0cm
\begin{Theo} {\rm (recursive definition of primitive hyperplanes and cocircuits of oriented cubes)}

Let  $\bf h$  be a primitive vector of  $C^{n-1}$  for which there exists  $b\in \N_0$   satisfying the following two conditions:
{\parindent=1cm 

(h-i)  $\forall  a\in \N_0, \ 0<a<|{\bf h}|,\ a\not=b \ \ S_a({\bf h}) \hookrightarrow S_b({\bf h})$

(h-ii) The signed set  $X_b:= ( \cup_{a<b} S_a({\bf h}),\ \cup_{a>b} S_a({\bf h}))$  is a  signed cocircuit of every oriented cube  $\mathcal M(C^{n-1})$}

Then, every vector ${\bf g}=({\bf h}, g)\in \N_0^{n}$, with  $0\leq g\leq b$, satisfies the following two properties:

{\parindent=1cm 
(g-i)  $\forall  a\in \N_0, \ 0<a<|{\bf g}|,\ a\not=b \ \ S_a({\bf g}) \hookrightarrow S_b({\bf g})$.

(g-ii) The signed set  $\tilde{X}_b= ( \cup_{a<b} S_a({\bf g}),\ \cup_{a>b} S_a({\bf g}))$  is a signed cocircuit of every oriented cube  $\mathcal M(C^{n})$.}
\end{Theo}

{\bf Proof.}} Let  $\tilde{\mathcal M}=\mathcal M(C^{n})$  be an oriented  cube. Consider vectors  ${\bf h}\in \N_0^{n-1}$  and   ${\bf g}=({\bf  h},g)\in \N_0^{n}$  in the conditions of the Theorem. Consider    $\tilde{H}:=cl(S_b({\bf g}))$. 

We prove the theorem in three steps. 

{\bf (I)}  $\tilde{\mathcal H}$  is a hyperplane of $\tilde{\mathcal M}$, whose pair of signed cocircuits restricts in the facet  $H_{n^+}$ of  $\tilde{\mathcal M}$  to the pair of signed cocircuit  $\pm X_b$  of  $S_b({\bf h})$.

{\bf (II)}  $\bf g$  satisfies  $(g-i)$. 

{\bf (III)} The unique extension of the cocircuit  $X_b$  of the facet  $H_{n^+}$ to a signed covector  of  $\tilde{\mathcal M}$  complementary of  $S_b({\bf g})$  and orthogonal to the rectangles of  $C^n$  is the signed set  $\tilde{X}_b$  defined  $(g-ii)$.\vskip 2mm

{\parindent=0cm
{\bf (I)} } Since  $\bf h$  satisfies  $(h-i)$  we know that  $S_{b-g}({\bf h})\hookrightarrow S_b({\bf h})$ implying also that  $S_{b-g}({\bf h})n^-\hookrightarrow S_b({\bf h})n^+$. Consequently  the ranks of  $S_b ({\bf g})$  and  $S_b ({\bf h})n^+ $  are related by:
$$rk (S_b ({\bf g}))= rk (S_b ({\bf h})n^+ )+1. \leqno{(1)}$$
The restriction of the canonically oriented cube  $\tilde{\mathcal M}$  to a facet is an oriented cube of rank  $rk(\tilde {\mathcal M})-1$.   The set  $H:= S_b ({\bf h})n^+$  is, by condition $(h-i)$  a hyperplane of the facet  $H_{n^+}$  of  $\tilde{\mathcal M}$  and therefore a hyperline of  $\tilde{\mathcal M}$  implying that   $rk(H)=rk(\tilde{\mathcal M})-2$. Replacing in  $(1)$ we have:
$$rk(S_b ({\bf g}))= rk(H)+1= rk(\tilde{\mathcal M})-1 \leqno{(2)}$$
Therefore  $\tilde{H} := cl_{\tilde{\mathcal M}}(S_b ({\bf g}))$  must be a hyperplane of  $\tilde{\mathcal M}$  and its complement the support of a pair  
 of signed cocircuits of  $\tilde{\mathcal M}$. Moreover, one of these signed cocircuits of  $\tilde{\mathcal M}$  must restrict in the facet  $H_{n^+}$  to the signed cocircuit  $X_b:= ( \cup_{a<b} S_a({\bf h})n^+,\ \cup_{a>b} S_a({\bf h})n^+)$ completing the proof of Step (I). \vskip 2mm
 
{\parindent=0cm
{\bf (II)} \it $\bf g$  satisfies  $(g-i)$.  }

Notice that for every  $a,\  0\leq a<g$,  $S_a({\bf g})=S_a({\bf h})n^+$ and since $\bf h$  satisfies  $(h-i)$  it is clear that in this case we have: $S_a({\bf g})n^+ \hookrightarrow S_b({\bf h})n^+$  and therefore $S_a({\bf g})n^+ \hookrightarrow S_b({\bf g})$.  For  $a>|{\bf h}|$,  since  $S_a({\bf g})=S_{a-g}({\bf h})n^-\hookrightarrow S_b({\bf h})n^+$  and the result in this case is also a direct consequence of  the fact that  $\bf h$  satisfies  $(h-i)$ 

For  $a,\ \ g\leq a \leq |{\bf h}|$, $S_a({\bf g})=S_a({\bf h})n^+\cup  S_{a-g}({\bf h})n^-$ intersects both facets  $H_{n^+}$  and  $H_{n^-}$. The hypothesis that  $\bf h$  satisfies  $(h-i)$  guarantees that $S_a({\bf g})n^+ \hookrightarrow S_b({\bf h})n^+$   and also that  $S_{a-g}({\bf h})n^-\hookrightarrow S_b({\bf h})n^+$  so, in order to conclude that   $S_a({\bf g}) \hookrightarrow S_b({\bf g})$  it is enough to prove  that there is a geometric realization of the 2-rectangle $r_2=(a=a; b=b)$  of the form  $R=( ({\bf u}_a,1)^+, ({\bf u}_{a-g},-1)^-, ({\bf v}_{b-g},-1)^-, ({\bf v}_b,1)^+)$   with  $({\bf u}_a,1)\in S_a({\bf h})n^+$,  $ ({\bf u}_{a-g},-1)\in S_{a-g}({\bf h})n^-$,  $({\bf v}_{b-g},-1)\in S_b({\bf h})n^-$   and  $ ({\bf v}_b,1))\in  S_b({\bf h})n^+$.

Consider the  $3$  or  $4$  rectangle of  $r'=(a-g<a; b-g<b)$  of  $\bf h$. The vector  $\bf h$  is a primitive vector of  $C^{n-1}$  so there is a geometric realization  $R'=({\bf u}_{a-g}^+, {\bf u}_a^-; {\bf v}_{b-g}^-, {\bf v}_b^+)$ of this rectangle in  $C^{n-1}$  leading the desired realization of the  2-rectangle  $r_2$  in  $C^n$  and concluding the proof that  ${\bf g}$  satisfies  $(g-i)$.\vskip 2mm

{\parindent=0cm
{\bf (III)} \it $\bf g$  satisfies  $(g-ii)$}

We know from  $(I)$  that there is a signed cocircuit  $\tilde{X}$  of   $\tilde{\mathcal M}$  complementary of the hyperplane $\tilde{H}$  whose restriction  to the facet  $H_{n^+}$  is  $X_b$.  The fact that  $\bf g$  satisfies  $(g-i)$, proved in (II),  implies directly  that the unique extension of  $X_b$  to  $H_{n^-}\setminus S_{b-g}({\bf h})n^-$  orthogonal to the rectangles of  $\tilde{\mathcal M}$  must satisfy  the conditions:  $\cup_{a<b}  S_a({\bf g}) \subseteq  \tilde{X}^+$  and  $\cup_{b< a\leq |\bf h|}  S_a({\bf g}) \subseteq  \tilde{X}^-$. 

Concerning the levels  $S_a({\bf g})$  with   $a>|{\bf h}|$  the fact that they are embedded in  $S_b({\bf g})$  implies that all the elements in each level will have the same sign however some more arguing is needed before concluding that they must all be negative in  $\tilde{X}$.

 We consider separately two cases:  {\it Case 1)}  $a\leq 2b$  and  {\it case 2)}  $a>2b$.\vskip 2mm

{\parindent=0cm
 {\it Case 1)}   In this case  consider the  3-rectangle of  $\bf g$: $r=(2b-a <b=b<a )$. If the rectangle  $r'=(b-g<a-g ; 2b-a<b)$  of  $\bf h$  is a  $3$  or  $4$ rectangle of  $\bf h$  then, since  $\bf h$  is primitive,  $r'$  is realizable.  Consider a realization  $R'= ({\bf u}_{b-g}, {\bf v}_{a-g}; {\bf v'}_{2b-a}, {\bf u'}_b)$  of   $r'$  in  $C^{n-1}$. Then clearly   $R= ({\bf v'}_{2b-a},1)^+\  ({\bf u'}_b,1)^-\ ({\bf u}_{b-g},-1)^-\ ( {\bf v}_{a-g}, -1))^+$  is a geometric realization  of   $r$  in  $C^n$  and by orthogonality with this circuit in any extension of  $X_b$  to $H_{n^-}\setminus S_{b-g}({\bf h})n^-$  the sign of  $( {\bf v}_{a-g}, -1))$  must be negative, implying that  $S_a(\bf g)\subseteq \tilde{X}^-$.  In the case  $r'$  is a $2$-rectangle of  $\bf h$, which occurs when  $g=a-b$, take  ${\bf u}_{b-g}\in S_{b-g}({\bf h})$  and  ${\bf v}_{b}\in S_{b}({\bf h})$. Clearly  $R=({\bf u}_{b-g},1)^+ \ ({\bf v}_b,1)^-\  ({\bf u}_{b-g},-1)^-\ ( {\bf v}_{b}, -1))^+$  is a geometric realization  of   $r$  in  $C^n$  and we also conclude that  $S_a(\bf g)\subseteq X^-$.\vskip 2mm

{\it Case 2)}} In this case  $r=(|{\bf h}|-(a-b) < b < |{\bf h}| < a )$  is a  $4$-rectangle of  $\bf g$ (notice that  $b<|{\bf h}|<a$)   and   $r'=(b-g<a-g; |{\bf h}|-(a-b)< |{\bf h}|)$  is a rectangle of   $\bf h$. If  $r'$  is a $3$ or $4$-rectangle of  $\bf h$  then, since  $\bf h$  is primitive, there is a realization  $({\bf u}_{b-g}, {\bf v}_{a-g}; {\bf v'}_{|{\bf h}|-(a-b)}, {\bf u'}_{|{\bf h}|})$  of   $r'$  in  $C^{n-1}$  and  $R= {\bf v'}_{|{\bf h}|-(a-b)},1)^+\  ({\bf u}_{b-g},-1)^- \ ({\bf u'}_{|{\bf h}|},1)^-\ ({\bf v}_{a-g}\  -1)^+)$  is geometric realization of  $r$. By orthogonality with this circuit we conclude as before that  $S_a({\bf g})\subseteq \tilde{X}^-$. The case  $r'$  is a 2-rectangle the argument is similar as in the previous case, leading to the conclusion that the unique extension  of  $X_b$  to  $H_{n^-}\setminus S_{b-g}({\bf h})n^-$  orthogonal to the rectangles of  $\tilde{\mathcal M}$ 
is the signed vector  $\tilde{X}_b$  of the theorem, thus proving that  $\bf g$  satisfies  $(g-ii)$.  $\square$}\vskip 5mm

The next Corollary whose proof is left to the reader restates Theorem 3.2 as a recursive procedure for constructing from signed cocircuits of every canonically oriented cube over $C^n$  signed cocircuits of every canonically oriented cube over  $C^{n+1}$. 

{\parindent=0cm  
\begin{Cor}
Let  $\bf h\in \N_0^n$  be a primitive vector such that the signed cocircuit of  $\mathcal Aff(C^n)$  defined by  $S_b({\bf h}):\  {\bf x}.{\bf h}= |{\bf h}|-2b$ 
is a signed cocircuit of every oriented cube $\mathcal M(C^n)$. For every  $c\leq b$ let  ${\bf g}:= ({\bf h},b-c)$, then the signed cocircuit  of  $\mathcal Aff(C^{n+1})$  defined by  $S_b({\bf g}):\   {\bf g}.{\bf x}= |{\bf g}|-2b$   is a signed cocircuit of every oriented cube  $\mathcal M(C^{n+1})$.
\end{Cor}
}\vskip 2mm

We are now ready to prove the main theorem.

\section{The main Theorem}
\vskip 2mm

\begin{Theo}{ For $n\leq 7$  the oriented matroid  $\mathcal Aff(C^n)$  is the unique oriented cube.}\end{Theo}

We recall  from \cite{dS2} that every oriented cube over  $C^n$  must have  $rank\ n+1$. The next Proposition about matroids then guarantees that in order to prove theorem 3.1 we "only" have to prove that every signed cocircuit of  $\mathcal Aff(C^n)$  must be a signed cocircuit of every oriented cube   $\mathcal M(C^n)$. 

{\parindent=0cm
\begin{Prop}
Consider two matroids  $\mathcal M$,   $\mathcal M'$, witout loops, over the same set  $E$ and with the same rank  $r$. Let  $\mathcal H,\mathcal H'$  denote the  families of hyperplanes  respectively of  $M$  and  $M'$  and assume that  $\mathcal H\subseteq \mathcal H'$.  Then  $\mathcal H =\mathcal H'$ and  $M=M'$.
\end{Prop}

{\bf Proof.}} Suppose that there is a hyperplane  $H'\in \mathcal H'\setminus \mathcal H$. Let  $B=\{h_1, \ldots, h_{r-1}\}$  be a basis of  $H'$  in  $M'$. Since $|B|<r$,  $B$  is contained in some hyperplane  $H\in \mathcal H$  of  $M$. The assumption that  $\mathcal H\subset \mathcal H'$  implies that  $H\in \mathcal H'$  and in  $\mathcal M'$  we have  $B\subset H\cap H'$. Now  $H\cap H'$  is a flat of  $M'$  whose rank is less then  the rank of  $H'$   contradicting the assumption that  $B$ is a basis of  $H'$. \vskip 5mm

The explicit definition of  $\mathcal  Aff(C^n)$  for  $n\leq 7$  that  we use is in terms of the family of signed cocircuits, defined by the following family  $\mathcal H_n$   of non-negative integer vectors, up to automorphisms of the class of orientations (i.e. of the unsigned arrangement of hyperplanes representing the oriented matroid)

{\parindent=0cm 
$\mathcal H_n:=\{(h_1\leq h_2\leq \ldots \leq h_n\leq h_{n+1})\in \N_0^{n+1}:\ \ gcd (h_1,\ldots, h_{n+1})=1\   and  $

$\ {\bf x}.(h_1,\ldots ,h_n)=h_{n+1}$ {\it defines\ a \ hyperplane \ of  \ $Aff(C^n) \}.$ }

Each vector  $(h_1\leq h_2\leq \ldots \leq h_n \leq h_{n+1})\in \mathcal H_n$  of the form  $(0^{\alpha_0}, h_1^{\alpha_1},\ldots \ h_k^{\alpha_k}) $  with  $\alpha_0\geq 0$  $\alpha_1,\ldots,\alpha_k\geq 1$  and   $\alpha_0 +\alpha_1+\ldots +\alpha_k= n+1$  determines exactly  $\frac{(n+1)!}{\alpha_0! \ldots \alpha_k!} 2^{n-\alpha_0}$  distinct hyperplanes (and pairs of signed cocircuits ) of  $\mathcal Aff(C^n)$.  

We recall that there is a natural procedure to generate vectors of  $\mathcal H_{n}$  from vectors of   $\mathcal H_{n-1}$  and  a recursive family $\mathcal G_n\subseteq \mathcal H_n$  that we briefly recall from    \cite{dS1}.}

{\parindent=0cm
\begin{Def} {\rm  (the family  $\mathcal G_n$) \cite{dS1}}

Given  ${\bf h}=(h_1\leq h_2\leq \ldots \leq h_n)\in \mathcal H_{n-1}$. For  $i=1,\ldots n$  let  ${\bf h_i}:= (h_1,\ldots, h_{i-1}, h_{i+1},\ldots , h_n)$, then the level  $S_b({\bf h_i})$  with  $b=	\frac{|{\bf h_i}|-h_i}{2}$  is a hyperplane of  $\mathcal Aff(C^{n-1})$  and for every  $c\leq b$, such that  $S_c({\bf h_i})\not=\emptyset$  the vector  ${\bf g}:= ({\bf h_i}, b-c, |{\bf h_i}|-b-c)$  is a vector of  $\mathcal H_n$. 

We denote by  $\mathcal G_n$  the family of vectors of  $\mathcal H_n$  which is obtained (up to reordering of entries) in this way.
\end{Def}}

The next Theorem gives the explicit definition of  $\mathcal H_n$  as well as the relation between  $\mathcal H_n$  and  $\mathcal G_n$, for  $n\leq 7$.\vskip 5mm

{\parindent=0cm
{\bf Theorem A}{(\cite{dS1})}

Consider  $\mathcal H_n:=\{ (h_1\leq \ldots \leq h_n\leq h_{n+1})\in \N_0^{n+1}:\  \ gcd (h_1,\ldots, h_{n+1})=1\   and\  (h_1,\ldots , h_n).{\bf x}=h_{n+1}\ defines\ a$  $\ hyperplane\ of\ \mathcal Aff(C^n)\}$.  The list  $\mathcal H_n$  as well as its sublist  $\mathcal G_n$  of definition 3.1  is the following   for  $n\leq  7$ :

$\mathcal H_1= \mathcal G_1= \{(1,1)\}$;

$\mathcal H_2=\mathcal G_2= \{(0,1,1)\}$;

$\mathcal H_3=\mathcal G_3= \{(0,0,1,1), (1,1,1,1)\}$;

$\mathcal H_4=\mathcal G_4= \{(0,0,0, 1,1), (0,1,1,1,1),(1,1,1,1,2)\}$;

$\mathcal H_5=\mathcal G_5= \{(0,0,0,0, 1,1), (0,0,1,1,1,1),(0,1,1,1,1,2), (1,1,1,1,1,1),$ 

\hskip 1.8cm $(1,1,1,1,1,3),(1,1,1,1,2,2), (1,1,1,2,2,3)\}$;

$\mathcal H_6=\mathcal G_6 =\{(0,0,0,0,0,1,1), (0,0,0,1,1,1,1),(0,1,1,1,1,2), (0,1,1,1,1,1,1),$

\hskip 8mm $(0,1,1,1,1,1,3),(0,1,1,1,1,2,2), (0,1,1,1,2,2,3), (1,1,1,1,1,1,2),$

\hskip 8mm $(1,1,1,1,1,1,4), (1,1,1,1,1,2,3), (1,1,1,1,2,2,2), (1,1,1,1,2,2,4),$

\hskip 8mm $(1,1,1,1,2,3,3), (1,1,1,1,3,3,4), (1,1,1,2,2,2,3), (1,1,1,2,2,2,5),$

\hskip 8mm  $(1,1,1,2,2,3,4), (1,1,1,2,3,3,5), (1,1,2,2,2,3,3), (1,1,2,2,3,3,4)$

\hskip 0.8cm $ (1,1,2,2,3,4,5)\}$;

$\mathcal H_7= \mathcal G_7 \cup$

\hskip 6mm $\{(1,1,2,4,4,5,6,7),(1,1,2,4,4,6,7,9), (1,2,3,3,4,4,5,6), (1,2,3,4,4,5,6,7),$

\hskip 6mm$(1,2,3,4,5,6,7,8), (1,3,3,3,4,4,5,5), (1,3,3,4,5,5,6,7), (2,2,2,2,3,3,3,5),$

\hskip 6mm $(2,2,2,3,3,3,4,7), (2,2,3,3,3,4,4,5), (2,2,3,3,4,4,5,7), (2,3,3,4,4,5,5,6)\}.$\vskip 1cm

\hskip 1cm We are now ready to prove the main Theorem. \vskip 5mm

{\bf Proof of Theorem  3.1.} }

Consider an oriented cube   $\mathcal M(C^n)$.  

In  order to prove that every signed cocircuit defined by a vector  $\bf h \in \mathcal H_n$  is a signed cocircuit of  $\mathcal M(C^n)$ it is enough to prove that one of the signed cocircuits  defined by an equation  ${\bf h_i.x}=h_i$  is a signed cocircuit of  $\mathcal M(C^n)$. This is a consequence of the behaviour of the families of signed sets, $\cal R$  and  $\cal F$, under the symmetries of the real cube (\cite {dS2}).

{\parindent=0cm
{\bf 1)} {\it Cases  $n\leq 6$.}

For  $n\leq 5$  and every  ${\bf h}\in \mathcal H_n$,  $\bf h_i$  is a primitive vector for every  $i=1,\ldots n+1$. Theorem 3.2. then guarantees that all the hyperplanes and corresponding signed cocircuits of  $\mathcal Aff(C^n)$  defined by vectors of   $\mathcal G_{n+1}= \mathcal H_{n+1}$  must be signed cocircuits of  $\mathcal M(C^n)$ and therefore, by Proposition 3.2  $\mathcal M(C^n)=\mathcal Aff(C^n)$.

{\bf 2)} {\it Case  $n= 7$.}

{\bf Claim 2.A)} {\it The signed cocircuits of  $\mathcal Aff(C^7)$  defined by the vectos of  $\mathcal G_7$  must be signed cocircuits of  $\mathcal M(C^7)$.}}

The arguments applied in the previous cases to vectors of  $\mathcal H_6$  lead to the conclusion that the hyperplanes and signed cocircuits of  $\mathcal Aff(C^7)$  defined by  $109$ of the $131$ vectors of  $\mathcal G_7$ must be hyperplanes and signed cocircuits of  $\mathcal M(C^7)$. The remaining  $22$  vectors of  $\mathcal G_7$ arise from the following  $5$  vectors of  $\mathcal H_6$: $(1,1,1,1,3,3,4),$  $  (1,1,2,2,2,3,3),$ $(1,1,2,2,3,3,4),$ $ (1,1,2,2,3,4,5)$  and  $(1,1,1,2,3,3,5).$

Each one of these  $5$ vectors has one restriction   ${\bf h_i}$  that is not primitive. In the case of the first  $4$  vectors, the restriction   $\bf h_1$ , in the case of the fifth the restriction  $\bf h_4$  so we can not apply  Theorem 3.2  to conclude that the hyperplanes of  $\mathcal G_n$   obtained  extending  the hyperplanes   $\bf h_i.x$ $=h_i=  S_b({\bf h_i)},\ with \ $ $b=\frac{|{\bf h_i}|-h_i}{2}$ in the facet  $H_{7^+}$  with an  $\bf h_i$-level in the facet   $H_{7^-}$  must be signed cocircuits  of  $\mathcal M(C^7)$ needs some further verification.

Actually in all the cases the proof goes in the same way, generalizing the proof of theorem 3.2:

First we verify that in all the  $5$  cases  the  $\bf h_i$-levels embed in the level  $S_b({\bf h_i}): \ {\bf h_i.x}= b$  ( $b=\frac{|{\bf h_i}|-h_i}{2}$). 

Secondly, for all  $c,\  0\leq c \leq b$  for which  ${\bf g}:= ({\bf h_i}, b-c, |{\bf h_i}|-b-c)$  is one of the  $109$ vectors of  $\mathcal G_n$   arising from primitive restrictions of vectors of  $\mathcal H_6$ there is nothing to prove. We retain those  $c's$ (e.g.  $c\geq h_7$)  for which  ${\bf g}$  is one of the remaining  $22$  vectors of  $\mathcal G_n$. For those we verify that $S_a({\bf g}) \hookrightarrow S_b({\bf g})$. Note that for that it is enough to prove that there is a  $2$-rectangle whose edges in each level connect the facets  $H_{7^+}$  and  $H_{7^-}$. 

The embeddings  $S_a({\bf g}) \hookrightarrow S_b({\bf g})$  imply: first, that  for every  ${\bf v}\in S_c({\bf h_i})7^-$ the hyperplane  $cl(S_b({\bf h_i})7^+ \cup {\bf  v})$  of  $\mathcal M(C^7)$   must contain all the elements of  $S_c({\bf h_i})7^-$.  Next, that the  only way of  extending the signed cocircuit  of  $S_b({\bf h_i})$  in the facet  $H_{7^+}$  to the complementary of $S_c({\bf h_i})7^-$  in the facet  $H_{7^-}$, orthogonally to the signed rectangles, must be the signed cocircuit of  $\mathcal Aff(C^7)$  complementary of  $S_b({\bf g})$. 

Once Claim 2A. is proved we start proving:\vskip 2mm

{\parindent=0cm
{\bf Claim 2.B)}{\it All the signed cocircuits of  $\mathcal Aff(C^7)$  determined by the  $12$  vectors of  $\mathcal H_7\setminus \mathcal G_7$  must be signed cocircuits of  $\mathcal M(C^7)$}}

Of these twelve vectors eight differ in exactly one entry from a vector of  $\mathcal G_7$, meaning that they span hyperplanes of  $\R^n$  parallel to the affine span of some hyperplane of  $\mathcal Aff(C^7)$  that we already know is a hyperplane of  $\mathcal M(C^7)$.

These eight vectors, together with the vector of  $\mathcal G_7$  differing in one entry are:

$(1,1,2,4,4,5,6,7), \ (1,2,3,3,4,4,5,6) \ and\ (1,2,3,4,4,5,6,7)$  all differing in one entry from  $(1,1,2,3,4,4,5,6)\in \mathcal G_7$,  

$(1,2,3,4,5,6,7,8)$  and  $(1,1,2,3,4,5,6,8)\in \mathcal G_7$,  

$(2,2,2,2,3,3,3,5)$  and  $(1,2,2,2,2,3,3,3)\in \mathcal G_7$,

$(2,2,2,3,3,3,4,7)$  and  $(1,2,2,2,3,3,3,4)\in \mathcal G_7$.

In all the cases, let  $\bf h$  be a restriction of the vector of  $\mathcal G_7$  defining a hyperplane  $S_c({\bf h})$  that we already know must be in  $\mathcal M(C^7)$  and let  $S_b({\bf h})$  be the level corresponding to the hyperplane of  $\mathcal Aff(C^7)$  whose pair of signed cocircuits we want to conclude must be in $\mathcal M(C^7)$.  In order to do so we proceed in all the cases in the same way:

{\parindent=0cm
(i) verify  that in  $\mathcal M$,  $rk(S_b({\bf h})) \leq 7$. For that, if necessary using the rectangles, we reduce to  $7$  the number of elements needed to span all the elements of that  $\bf h$-level.

(ii) verify that for all $a\not=b$, $S_a({\bf h }) \hookrightarrow S_b({\bf h})$. This implies, in particular, that the hyperplane  $S_c({\bf h})$  is embedded in  $S_b({\bf h})$  and consequently, for an element  ${\bf v}_a\in S_a$  we have  $rk(cl (S_b \cup {\bf v}_a))\leq rk(S_b) +1$.  On the other hand the embedding $S_a\hookrightarrow S_b$  implies that  $rk(cl (S_b \cup {\bf v}_a))= rk(cl (S_b \cup S_a))= rk(\mathcal M)=8 $  and we conclude that  $rk(S_b)=7$  and  $cl(S_b)$  must be a hyperplane of  $\mathcal M$.

(iii) the fact that for all $a\not=b$, $S_a \hookrightarrow S_b$ implies that in the signed cocircuit complementary of 
$cl(S_b)$ all the elements of  $S_a$ must all have the same sign, or be all contained in the hyperplane $cl(S_b)$. In the eight cases a "ladder" of signed rectangles leads the conclusion that the pair of signed cocircuits of  $\mathcal Aff(C^7)$  complementary of  $S_b$  is the unique signature of the complement of  $S_b$  orthogonal to the rectangles of  $C^7$  and therefore, these cocircuits of  $\mathcal Aff(C^7)$  must be cocircuits of  $\mathcal M$.}

Of the four remaining vectors of  $\mathcal H_7\setminus \mathcal G_7$  the following two:  

$(1,1,2,4,4,6,7,9)$  and  
$(1,2,3,3,4,5,5,6)$  differ in exactly one entry from the vectors, respectively, $(1,1,2,4,4,5,6,7)$  and  $(1,2,3,3,4,4,5,6)$ of  $\mathcal H_7$  whose signed cocircuits we already know must be in  $\mathcal M$ and we proceed exactly as before in order to conclude that the corresponding signed cocircuits of  $\mathcal Aff(C^7)$  must be signed cocircuits of  $\mathcal M$ too.

Finally we are left with the two vectors  $(1,3,3,3,4,4,5,5)$  and  

\hskip -1cm $(1,3,3,4,5,5,6,7)$. In both cases we prove that the pair of signed cocircuits of the hyperplane defined by  $S_b({\bf h_8}):=\ \ {\bf h_8.x}=h_8$  ($b=\frac{|{\bf h_8}|-h_8}{2}$)  must be cocircuits of  $\mathcal M$. In both cases the hyperplane contains exactly  seven elements of  $C^7$, therefore we know that  $rk(cl( S_b))\leq 7$.

To conclude that  $cl(S_b)=S_b$  is a hyperplane whose complement, signed as in  $\mathcal Aff(C^7)$  must be signed cocircuits of  $\mathcal M$  we proceed as before: first, we verify that  $\forall a\not=b,\ S_a\hookrightarrow S_b$. Next use the net of signed rectangles to conclude that if one element of a different level  $S_a$   belongs to  $cl(S_b)$  then  $cl(S_b)=C^7$, implying that $cl(S_b)=S_b$. Then, essentially the same ladder of rectangles proves the unicity of the signature, concluding the Proof of Theorem 3.1.

\section{Final Remarks}

In this note we defined recursively a family of nonnegative integer vectors orthogonal to hyperplanes of the real cube - primitive vectors. We proved that the corresponding primitive  hyperplanes and signed cocircuits of  the real cube  $\mathcal Aff(C^n)$  must be hyperplanes and signed cocircuits of every oriented cube.  The proof relied upon a direct characterization of the behaviour of these nonnegative integer vectors with respect to the net of rectangles of  $\mathcal Aff(C^n)$. 

Although explicit descriptions of the affine cube for  $n=8$  have been computed  \cite{AA}, and our methods apply to this case,  
an approach of this next case without computer aid would still be too lengthy. A better understanding of the real affine cube and its symmetries is needed.

For instance, note that the hyperplanes defined by primitive vectors are strictly contained in the family  ${\mathcal G}_n$ (definition 3.1) which for  $n\geq 7$  is strictly contained in the family  $\mathcal H_n$. Although we think it must be true, we were not able to prove that the  hyperplanes and signed cocircuits determined by the vectors of  $\mathcal G_n$  must be hyperplanes and signed cocircuits of every oriented cube. This would, in particular,  simplify further our proof of the case  $n=7$.

{\parindent=0cm

{\bf E. Gioan} 

{e-mail adress:} emeric.gioan@lirmm.fr

LIRMM (UMR 5506-CC477), 161,r. Ada 34095 Montpellier Cedex 5 - FRANCE\vskip 5mm

{\bf I. P. Silva} 

{e-mail adress:} ipsilva@fc.ul.pt

Dep. de Matem\'atica, Faculdade de Ciencias da Universidade de Lisboa

Edicio C6, Campo Grande, 1749-016 Lisboa, PORTUGAL }

\end{document}